\pdfoutput=1

\documentclass[reqno]{amsart}

\usepackage[T1]{fontenc}
\usepackage{lmodern}

\usepackage[utf8]{inputenc}
\usepackage[english]{babel}
\usepackage{microtype} 

\usepackage{amsmath,amssymb,amsthm,mathrsfs,latexsym,mathtools,mathdots,booktabs,enumerate,tikz,bm,comment,tikz,hyperref}
\usepackage[centertableaux]{ytableau}

\usetikzlibrary{arrows,positioning,shapes.geometric}

\definecolor{pLightBlue}{RGB}{150,200,255}
\definecolor{pLightOrange}{RGB}{231,183,100}
\definecolor{pBlue}{RGB}{86,139,190}
\definecolor{pCyan}{RGB}{149,186,201}
\definecolor{pSand}{RGB}{184,166,121}
\definecolor{pAlgae}{RGB}{87,115,135}
\definecolor{pSkin}{RGB}{236,216,167}
\definecolor{pGray}{RGB}{156,175,156}
\definecolor{pPink}{RGB}{215,114,127}
\definecolor{pOrange}{RGB}{211,153,80}

\usepackage{graphicx}
\DeclareGraphicsExtensions{.pdf}

\newtheorem{theorem}{Theorem}
\newtheorem{proposition}[theorem]{Proposition}
\newtheorem{lemma}[theorem]{Lemma}

\theoremstyle{definition}
\newtheorem{definition}[theorem]{Definition}
\newtheorem{example}[theorem]{Example}
\newtheorem{remark}[theorem]{Remark}

\newtheorem{problem}[theorem]{Problem}

\newcommand*\circled[1]{\tikz[baseline=(char.base)]{\node[shape=circle,draw,inner sep=0.1pt,minimum size=1.3em] (char) {#1};}}

\newcommand*\posetNode[2]{$\circled{$#1$}_{#2}$}

\newcommand{\sgn}{\operatorname{sgn}}

\newcommand{\defin}[1]{%
\relax\ifmmode%
\textcolor{blue}{#1}%
\else \textcolor{blue}{\emph{#1}}%
\fi%
}

\newcommand{\setN}{\mathbb{N}}
\newcommand{\setP}{\mathbb{N}^{+}}

\newcommand{\xvec}{\mathbf{x}}

\newcommand{\tagg}{\mathfrak{t}}

\newcommand{\opsurj}{\mathcal{O}}

\newcommand{\poseta}{a}
\newcommand{\posetb}{b}

\newcommand{\powerSum}{\mathrm{p}}
\newcommand{\schur}{\mathrm{s}}
\newcommand{\gessel}{\mathrm{F}}
\newcommand{\psih}{\hat{\psi}}

\DeclareMathOperator{\wt}{wt}
\DeclareMathOperator{\length}{\ell}

\DeclareMathOperator{\rot}{root}

\newcommand{\thsup}{\textnormal{th}}

\tikzset{every picture/.append
  style={scale=1,
	baseline=(current bounding box.center),
	x=1em,
	y=1em,
	naturalEdge/.style={line width=0.8pt},
	strictEdge/.style={line width=0.7pt,line join=round,double=white},
	entries/.style={xshift=-0.5em,yshift=-0.5em,font=\small},
	circled/.style={circle,draw,inner sep=0pt,minimum size=1em},
	bigcircled/.style={circle,draw,inner sep=0pt,minimum size=2.5em},
	hookTreeNode/.style={circle,draw,inner sep=0.1pt,minimum size=3.3em},
	hookTreeNodeInternal/.style={circle,draw,fill=pLightBlue,inner sep=0.1pt,minimum size=3.3em}
	}
}

\author{Per Alexandersson}
\email{per.w.alexandersson@gmail.com}
\address{Department of Mathematics,
Stockholm University,
S-10691, Stockholm, Sweden}

\author{Olivia Nabawanda}
\email{nabawandaolivia100@gmail.com}
\address{Department of Mathematics,
Mbarara University of Science \& Technology, Mbarara, Uganda}


\title{A weighted Murnaghan--Nakayama rule for $(P, \omega)$-partitions}

\begin{document}

\begin{abstract}
 The $(P, \omega)$-partition generating function $K_{(P, \omega)}(\xvec)$ is a quasisymmetric function obtained from a labeled poset.
 Recently, Liu and Weselcouch gave a formula for the coefficients of $K_{(P, \omega)}(\xvec)$
 when expanded in the quasisymmetric power sum function basis.
 This formula generalizes the classical Murnaghan--Nakayama rule for Schur functions.
 
 We extend this result to \emph{weighted} $(P, \omega)$-partitions and provide a short combinatorial proof,
 avoiding the Hopf algebra machinery used by Liu--Weselcouch.
\end{abstract}

\maketitle

\section{Introduction}

Given a finite poset $P$ with a labeling $w: P\rightarrow \{1,2,\dotsc,n\}$, the pair $(P, \omega)$ is referred
to as a \defin{labeled poset}.
For each $(P, \omega)$, there is a quasisymmetric
function $K_{(P, \omega)}(\xvec)$ which enumerates 
certain order-preserving maps $f: P\rightarrow \mathbb{N}^{+}$.
There are several natural bases for the space of quasisymmetric functions, many of which
generalize bases for the space of symmetric functions.
One such basis is the quasisymmetric power sums introduced by C.~Ballantine~et~al.~\cite{BallantineDaughertyHicksMason2020}.

The expansion of $K_{(P, \omega)}(\xvec)$ for naturally labeled posets
was done by P.~Alexandersson and R.~Sulzgruber in \cite{AlexanderssonSulzgruber2019} (Theorem~\ref{thm:ASFormula} below),
and the expansion for \emph{arbitrary labeled posets} was given 
in R.~Liu and M.~Weselcouch~\cite{LiuWeselcouch2020} (Theorem~\ref{thm:liu} below).
In particular, the classical Murnaghan--Nakayama rule for the expansion of
Schur functions into power sum symmetric functions is a special case of the Liu--Weselcouch formula; see Section~\ref{sec:MurnaghanNakayamaConnection}.

The quasisymmetric power sum expansion of a \emph{weighted generalization} of $K_{(P, \omega)}(\xvec)$ for naturally labeled
posets was also given in \cite{AlexanderssonSulzgruber2019} (Theorem~\ref{thm:ASFormula2} below).
Weighted $P$-partitions are a natural generalization of $P$-partitions.
We let $K^d_{(P, \omega)}(\xvec)$ denote the weighted generalization of $K_{(P, \omega)}(\xvec)$,
where $d:P \to \setP$ is a weight on the vertices of $P$.
For example, in \cite{AliniaeifardWangWilligenburg2023}, the authors also consider weighted $P$-partitions,
and use this as motivation to introduce and study what they 
call the \emph{combinatorial quasisymmetric power sum functions} (independently discovered by A.~Lazzeroni~\cite{Lazzeroni2023}). As another example, D.~Grinberg and E.~Vassilieva introduce a basis for the space of quasisymmetric
functions using weighted $P$-partitions in \cite{GrinbergVassilieva2021}.
Weighted colorings in general have also been studied extensively; see for example \cite{CrewSpirkl2020,AlistePrietoCrewSpirklZamora2021,AliniaeifardWangvanWilligenburg2021}.
\medskip

The main result in this paper is a unifying formula (Theorem~\ref{thm:main} below)
which provides the quasisymmetric power sum expansion of
 $K^d_{(P, \omega)}(\xvec)$, thus providing a Murnaghan--Nakayama rule for weighted labeled posets.
In particular, we give a new short proof of the quasisymmetric power sum expansions given in
\cite{AlexanderssonSulzgruber2019} and \cite{LiuWeselcouch2020}.
Our method is purely combinatorial and does not utilize Hopf algebras or the fundamental quasisymmetric functions.

\medskip
The paper is outlined as follows.
In Section~\ref{sec:prelim} we give some preliminary definitions and properties of the
quasisymmetric functions $K^d_{(P, \omega)}(\xvec)$. The main properties are
two recurrence relations satisfied by the $K^d_{(P, \omega)}(\xvec)$'s.

In Section~\ref{sec:mainresult}, the quasisymmetric power sum expansion of $K^d_{(P,\omega)}(\xvec)$
is given. With the help of the two aforementioned recurrence relations, we
use an inductive argument to reduce to the case of naturally labeled weighted chains.
We give a new proof in this case using a probabilistic approach.
Interestingly, this has a connection with linear extensions of certain trees.

Finally in Section~\ref{sec:problems}, we discuss some remaining open problems.

\section{Preliminaries}\label{sec:prelim}

\subsection{Labeled posets and order-preserving maps}

Let $P$ be a poset (partially ordered set) on $n$ elements.
A \defin{labeling} is a bijection $\omega : P \to \{1, 2, \dotsc, n\}$ 
and the pair $(P, \omega)$ is called a \defin{labeled poset}.
We say that an edge $(a, b)$ with $a<_P b$ is \defin{strict} if $\omega(a) > \omega(b)$. 
Otherwise, we say that an edge is \defin{weak} or \defin{natural}. 
In all our illustrations, we represent strict edges with double lines while
weak edges are represented as regular lines. 
In all our applications, only the edges in the Hasse diagram of $P$
matter and we mainly focus on these.

If all the edges of a labeled poset $(P, \omega)$ are natural, then the poset is said to be \defin{naturally labeled}. 
Note that every poset $P$ admits (at least) one \defin{natural labeling} $\omega$.

For example, the labeled poset in \eqref{eq:examplePoset}
has the strict edges $(7,1)$, $(7,2)$, $(3,2)$, $(6,4)$, $(9,8)$, while the remaining edges are weak.
	\begin{equation}\label{eq:examplePoset}
	(P, \omega) =
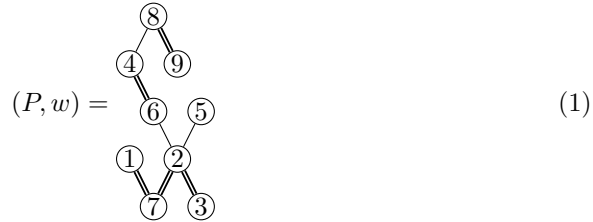
	\begin{tikzpicture}[scale=0.9]
	\draw
	(1,1)node[circled](v1){7}
	(3,1)node[circled](v2){3}
	(2,3)node[circled](v3){2}
	(1,5)node[circled](v4){6}
	(3,5)node[circled](v5){5}
	(0,7)node[circled](v6){4}
	(0,3)node[circled](v7){1}
	(2,7)node[circled](v8){9}
	(1,9)node[circled](v9){8}
	;
	\draw (v7)--(v1)--(v3)--(v4)--(v6)--(v9);
	\draw (v3)--(v5);
	\draw[strictEdge] (v1)--(v3);
	\draw[strictEdge] (v2)--(v3);
	\draw[strictEdge] (v4)--(v6);
	\draw[strictEdge] (v1)--(v7);
	\draw[strictEdge] (v8)--(v9);
	\end{tikzpicture}
	\end{equation} 

\begin{definition}
Let $(P,\omega)$ be a labeled poset. A \defin{$(P, \omega)$-partition} 
$f$ is a map $f: P \to \setP$ satisfying the following two properties:
\begin{enumerate}
 \item $\poseta <_P \posetb \implies f(\poseta) \leq f(\posetb)$  and
 \item $\poseta <_P \posetb \text{ and } \omega(\poseta) > \omega(\posetb) \implies f(\poseta) < f(\posetb)$.
\end{enumerate}
If $P$ is naturally labeled, then we simply call $f$ a \defin{$P$-partition}.
\end{definition}

\begin{example}\label{ex:a}
The following maps $f_1$ and $f_2$ are examples of $(P, \omega)$-partitions 
for the labeled poset in \eqref{eq:examplePoset},
where $f(\poseta)$ is written at vertex $\poseta \in P$:
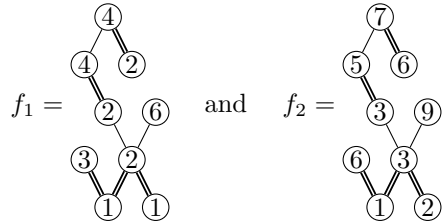
\begin{equation*}
f_1 = 
\begin{tikzpicture}[scale=0.9]
\draw
(1,1)node[circled](v1){1}
(3,1)node[circled](v2){1}
(2,3)node[circled](v3){2}
(1,5)node[circled](v4){2}
(3,5)node[circled](v5){6}
(0,7)node[circled](v6){4}
(0,3)node[circled](v7){3}
(2,7)node[circled](v8){2}
(1,9)node[circled](v9){4}
;
\draw (v7)--(v1)--(v3)--(v4)--(v6)--(v9);
\draw (v3)--(v5);
\draw[strictEdge] (v1)--(v3);
\draw[strictEdge] (v2)--(v3);
\draw[strictEdge] (v4)--(v6);
\draw[strictEdge] (v1)--(v7);
\draw[strictEdge] (v8)--(v9);
\end{tikzpicture}
\quad \text{ and } \quad
f_2 = 
\begin{tikzpicture}[scale=0.9]
\draw
(1,1)node[circled](v1){1}
(3,1)node[circled](v2){2}
(2,3)node[circled](v3){3}
(1,5)node[circled](v4){3}
(3,5)node[circled](v5){9}
(0,7)node[circled](v6){5}
(0,3)node[circled](v7){6}
(2,7)node[circled](v8){6}
(1,9)node[circled](v9){7}
;
\draw (v7)--(v1)--(v3)--(v4)--(v6)--(v9);
\draw (v3)--(v5);
\draw[strictEdge] (v1)--(v3);
\draw[strictEdge] (v2)--(v3);
\draw[strictEdge] (v4)--(v6);
\draw[strictEdge] (v1)--(v7);
\draw[strictEdge] (v8)--(v9);
\end{tikzpicture}
\end{equation*}
\end{example}

\subsection{Quasisymmetric functions}

\begin{definition}
A \defin{quasisymmetric function} in variables $\xvec = (x_1, x_2, \ldots)$ is a formal 
power series $g = g(\xvec) \in \mathbb{Q}[[x_1, x_2, \ldots]]$ of finite degree such 
that for every sequence $d_1, d_2, \dotsc, d_k$ of positive integer exponents, 
the coefficient of $x_{1}^{d_1}x_{2}^{d_2}\cdots x_{k}^{d_k}$ in $g$
is equal to the coefficient of $x_{i_1}^{d_1}x_{i_2}^{d_2}\cdots x_{i_k}^{d_k}$
whenever $i_1< i_2<\cdots < i_k$.
\end{definition}

Basis elements in the space of homogeneous quasisymmetric functions (of degree $n$)
are usually indexed by integer compositions $\alpha = (\alpha_1,\dotsc,\alpha_\ell)$,
where each $\alpha_i$ is a positive integer (and $\alpha_1+\dotsb+\alpha_\ell=n$).
One such basis is the \defin{monomial quasisymmetric functions}, defined as 
\begin{equation}\label{eq:qmonomDef}
 \defin{M_\alpha(\xvec)} \coloneqq \sum_{i_1 < i_2 < \dotsb < i_\ell} 
 x_{i_1}^{\alpha_1} x_{i_2}^{\alpha_2} \dotsm x_{i_\ell}^{\alpha_\ell}
\end{equation}
where $\alpha = (\alpha_1,\dotsc,\alpha_\ell)$.

\smallskip 

The generating function of all $(P, \omega)$-partitions $f$ is defined as
\[
\defin{K_{(P, \omega)}(\xvec)} \coloneqq  \sum_{f} \prod_{a\in P} x_{f(a)}.
\] 
For example, the functions $f_1$ and $f_2$ in Example~\ref{ex:a} 
contribute $x_1^2 x_2^3 x_3 x_4^2 x_6$ and $x_1x_2x_3^2x_5x_6^2x_7x_9$ to the 
sum for computing $K_{(P, \omega)}(\xvec)$.

It is now immediate from the definitions that $K_{(P, \omega)}(\xvec)$ 
is a (homogeneous) quasisymmetric function of degree $n$, where $n$ is the number of elements in $P$.
Quasisymmetric functions were first studied explicitly by I.~Gessel~\cite{Gessel1984}, 
who drew the inspiration to study them from the generating function $K_{(P, \omega)}(\xvec)$. 
Since then, quasisymmetric functions have become a prominent 
area of research in algebraic combinatorics. 
For more background, applications and recent results on quasisymmetric functions; see 
\cite[Section 7.19]{StanleyEC2} or \cite{LuotoEtAl2013IntroQSymSchur,Mason2019}.

\subsection{Weighted \texorpdfstring{$(P,\omega)$}{(P,omega)}-partitions}

A \defin{weighted poset} $(P, \omega, d)$ consists of a labeled poset $(P, \omega)$ 
and assignment of positive integer \defin{weights} $d:P \to \setP$ to the vertices of $P$.
We can now consider a weighted generalization of $K_{(P, \omega)}(\xvec)$.
We let \defin{$K_{(P, \omega)}^d(\xvec)$} be the generating function defined as
\[
\defin{ K_{(P, \omega)}^d(\xvec) } \coloneqq \sum_{f} \prod_{a\in P} x_{f(a)}^{d(a)}
\]
where the sum ranges over all $(P, \omega)$-partitions.

Note that this is a homogeneous quasisymmetric function of degree $\sum_{a \in P} d(a)$.
\begin{example}
Consider the labeled poset $P$ in Example~\ref{ex:a} in which each vertex $a\in P$ 
has been assigned a weight, $d(a)$. Below is a $(P, \omega)$-partition $f$,
and it contributes $x_{1}^{2}x_{2}^{1}x_{3}^{3+5}x_{4}^{8+1}x_{5}^{6+4+2}$ 
to the sum $K_{(P, \omega)}^d(\xvec)$:
\[ 
\begin{tikzpicture}[scale=1.7]
\draw
(2,9)node(v9){\posetNode{5}{d(\cdot)=2}}
(4,7)node(v8){\posetNode{3}{d(\cdot)=5}}
(-2,7)node(v7){\posetNode{5}{d(\cdot)=4}}
(-2,3)node(v6){\posetNode{4}{d(\cdot)=1}}
(2,1)node(v1){\posetNode{2}{d(\cdot)=1}}
(6,1)node(v2){\posetNode{1}{d(\cdot)=2}}
(4,3)node(v3){\posetNode{3}{d(\cdot)=3}}
(2,5)node(v4){\posetNode{4}{d(\cdot)=8}}
(6,5)node(v5){\posetNode{5}{d(\cdot)=6}};
\draw[strictEdge] (v1)--(v3);
\draw[strictEdge] (v2)--(v3);
\draw (v3)--(v4);
\draw (v3)--(v5);
\draw[strictEdge] (v1)--(v6);
\draw[strictEdge] (v4)--(v7);
\draw[strictEdge] (v8)--(v9);
\draw (v7)--(v9);
\end{tikzpicture}
\]
Here, $d(\cdot)=6$ indicates that the vertex has weight $6$.
\end{example}

The generating functions $K_{(P, \omega)}^d(\xvec)$ satisfy several recursions,
as we shall see in the following two lemmas.
\begin{lemma}\label{lem:addEdge}
Suppose $(P, \omega,d)$ is a weighted poset where $a,b \in P$ are incomparable and $\omega(a)<\omega(b)$.
Let $P'$ and $P''$ be obtained from $P$ by adding the relation $a <_{P'} b$ and $a>_{P''} b$, respectively.
Then the generating function $K_{(P,\, \omega)}^d(\xvec)$
satisfies the recurrence relation
\[
K_{(P,\,\omega)}^d(\xvec) = K_{(P',\, \omega)}^{d}(\xvec) + K_{(P'',\, \omega)}^{d}(\xvec).
\]
\end{lemma}
\begin{proof}
Any $(P, \omega)$-partition contributing to $f$ either has $f(a) \leq f(b)$ or $f(a)>f(b)$.
By breaking up $K_{(P,\omega)}^d(\xvec)$ into these two cases, we get the two generating functions on the right.
Pictorially, we can illustrate this as below
(dashed edges illustrate relations only involving one of $a$ and $b$): 
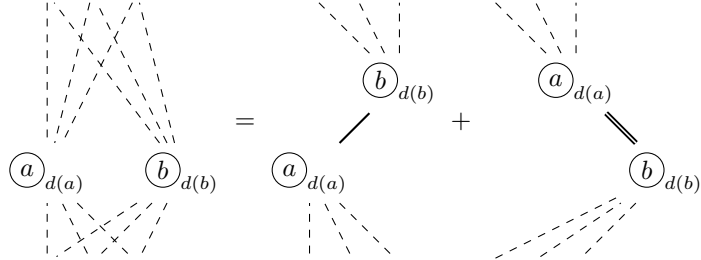
\[
\begin{tikzpicture}[scale=1.7]
\draw
(-2,0) node(extremeBL){}
( 4,0) node(extremeBR){}
(-2,6) node(extremeTL){}
( 4,6) node(extremeTR){}
(0,2) node(v1){\posetNode{a}{d(a)}}
(2,2) node(v2){\posetNode{b}{d(b)}}
(0,6)node(topLeft){ }
(1,6)node(topMid){ }
(2,6)node(topRight){ }
(0,0)node(botLeft){ }
(1,0)node(botMid){ }
(2,0)node(botRight){ };
\draw[gray,thick] (botLeft)--(v1);
\draw[gray,thick] (botMid)--(v1);
\draw[gray,thick] (botRight)--(v1);
\draw[gray,thick] (topLeft)--(v1);
\draw[gray,thick] (topMid)--(v1);
\draw[gray,thick] (topRight)--(v1);
\draw[gray,thick] (v2)--(topLeft);
\draw[gray,thick] (v2)--(topMid);
\draw[gray,thick] (v2)--(topRight);
\draw[gray,thick] (v2)--(botLeft);
\draw[gray,thick] (v2)--(botMid);
\draw[gray,thick] (v2)--(botRight);
\draw[dashed,blue] (v1)--(extremeBL);
\draw[dashed,blue] (v1)--(extremeTL);
\draw[dashed,red] (v2)--(extremeBR);
\draw[dashed,red] (v2)--(extremeTR);
\draw (1,-0.5) node {$\in (P,\omega)$};
\end{tikzpicture}
=
\begin{tikzpicture}[scale=1.7]
\draw
(-2,0) node(extremeBL){}
( 4,0) node(extremeBR){}
(-2,6) node(extremeTL){}
( 4,6) node(extremeTR){}
(0,2) node(v1){\posetNode{a}{d(a)}}
(2,4) node(v2){\posetNode{b}{d(b)}}
(0,6)node(topLeft){ }
(1,6)node(topMid){ }
(2,6)node(topRight){ }
(0,0)node(botLeft){ }
(1,0)node(botMid){ }
(2,0)node(botRight){ };
\draw[naturalEdge] (v1)--(v2);
\draw[gray,thick] (botLeft)--(v1);
\draw[gray,thick] (botMid)--(v1);
\draw[gray,thick] (botRight)--(v1);
\draw[gray,thick] (v2)--(topLeft);
\draw[gray,thick] (v2)--(topMid);
\draw[gray,thick] (v2)--(topRight);
\draw[dashed,blue] (v1)--(extremeBL);
\draw[dashed,blue] (v1)--(extremeTL);
\draw[dashed,red] (v2)--(extremeBR);
\draw[dashed,red] (v2)--(extremeTR);
\draw (1,-0.5) node {$\in (P',\omega)$};
\end{tikzpicture}
+
\begin{tikzpicture}[scale=1.7]
\draw
(-2,0) node(extremeBL){}
( 4,0) node(extremeBR){}
(-2,6) node(extremeTL){}
( 4,6) node(extremeTR){}
(0,4) node(v1){\posetNode{a}{d(a)}}
(2,2) node(v2){\posetNode{b}{d(b)}}
(0,6)node(topLeft){ }
(1,6)node(topMid){ }
(2,6)node(topRight){ }
(0,0)node(botLeft){ }
(1,0)node(botMid){ }
(2,0)node(botRight){ };
\draw[strictEdge] (v1)--(v2);
\draw[gray,thick] (botLeft)--(v2);
\draw[gray,thick] (botMid)--(v2);
\draw[gray,thick] (botRight)--(v2);
\draw[gray,thick] (v1)--(topLeft);
\draw[gray,thick] (v1)--(topMid);
\draw[gray,thick] (v1)--(topRight);
\draw[dashed,blue] (v1)--(extremeBL);
\draw[dashed,blue] (v1)--(extremeTL);
\draw[dashed,red] (v2)--(extremeBR);
\draw[dashed,red] (v2)--(extremeTR);
\draw (1,-0.5) node {$\in (P'',\omega)$};
\end{tikzpicture}
\]
\end{proof}

\begin{lemma}\label{lem:splitWeight}
Let $(P, \omega, d)$ be a weighted poset and $a \in P$ be any
vertex for which $d(a) = d_1 + d_2$
for some positive integers $d_1$ and $d_2$.
The generating function $K_{(P,\, \omega)}^d(\xvec)$ satisfies the relation
\[
K_{(P,\, \omega)}^d(\xvec) = K_{(P',\, \omega')}^{d'}(\xvec) - K_{(P',\, \omega'')}^{d'}(\xvec),
\] 
where we obtain $P'$ from $P$ by splitting vertex $a$ into two new vertices $a_- <_{P'} a_+$, 
each inheriting all relations from $P$. 
If $\omega(a)=i$ in $P$, then we obtain two new labelings $\omega'$ and $\omega''$
using the labels $1<2<\dotsb < i-1 < i_- < i_+ < i+1 < \dotsb < n$,
with
\begin{itemize}
 \item $\omega'(a_-)=i_-$, $\omega'(a_+)=i_+$ and $\omega'$ agreeing with $\omega$ on all other vertices; and
 \item $\omega''(a_-)=i_+$, $\omega''(a_+)=i_-$ and $\omega''$ agreeing with $\omega$ on all other vertices.
\end{itemize}
Finally, we let $d'$ be defined so that $d(a_-)=d_1$, $d(a_+)=d_2$
and $d'$ agrees with $d$ outside $a_-$ and $a_+$;
see Figure~\ref{fig:splitWeight}.
\end{lemma}
\begin{proof}
The $(P, \omega)$-partitions contributing to the first term in the right-hand side 
satisfy $f(a_-) \leq f(a_+)$. 
The cases where we have $f(a_-) < f(a_+)$ are 
counted by the second term in the right-hand side.
The difference is therefore the cases where $f(a_-) = f(a_+)$
and this is exactly the left-hand side; see Figure~\ref{fig:splitWeight}.
\begin{figure}[!ht]
 \centering
 \begin{tikzpicture}[scale=1.7]
\draw
(1,3)node(v1){\posetNode{i}{d_1+d_2}}
(0,6)node(topLeft){ }
(1,6)node(topMid){ }
(2,6)node(topRight){ }
(0,0)node(botLeft){ }
(1,0)node(botMid){ }
(2,0)node(botRight){ };
\draw[gray,thick] (botLeft)--(v1)--(topLeft);
\draw[gray,thick] (botMid)--(v1)--(topMid);
\draw[gray,thick] (botRight)--(v1)--(topRight);
\end{tikzpicture}
=
 \begin{tikzpicture}[scale=1.7]
\draw
(0,2) node(v1){\posetNode{i_-}{d_1}}
(2,4) node(v2){\posetNode{i_+}{d_2}}
(0,6)node(topLeft){ }
(1,6)node(topMid){ }
(2,6)node(topRight){ }
(0,0)node(botLeft){ }
(1,0)node(botMid){ }
(2,0)node(botRight){ };
\draw[naturalEdge] (v1)--(v2);
\draw[gray,thick] (botLeft)--(v1);
\draw[gray,thick] (botMid)--(v1);
\draw[gray,thick] (botRight)--(v1);
\draw[gray,thick] (v2)--(topLeft);
\draw[gray,thick] (v2)--(topMid);
\draw[gray,thick] (v2)--(topRight);
\end{tikzpicture}
$-$
 \begin{tikzpicture}[scale=1.7]
\draw
(2,4) node(v1){\posetNode{i_-}{d_2}}
(0,2) node(v2){\posetNode{i_+}{d_1}}
(0,6)node(topLeft){ }
(1,6)node(topMid){ }
(2,6)node(topRight){ }
(0,0)node(botLeft){ }
(1,0)node(botMid){ }
(2,0)node(botRight){ };
\draw[strictEdge] (v1)--(v2);
\draw[gray,thick] (botLeft)--(v2);
\draw[gray,thick] (botMid)--(v2);
\draw[gray,thick] (botRight)--(v2);
\draw[gray,thick] (v1)--(topLeft);
\draw[gray,thick] (v1)--(topMid);
\draw[gray,thick] (v1)--(topRight);
\end{tikzpicture}
\caption{Splitting a vertex with weight $d_1+d_2$ into two smaller vertices.}\label{fig:splitWeight}
\end{figure}
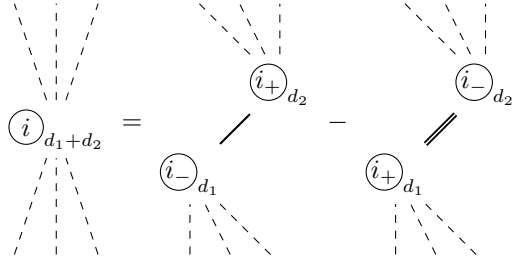
Note that all other strict and weak inequalities in the labeled posets are preserved,
since we do not change the relative order of any labels involving at most one of $a_-$ and $a_+$.
\end{proof}

\begin{remark}\label{rem:recursion}
Note that we may use the two above relations to 
express any $K_{(P, \omega)}^d(\xvec)$ as linear combinations of such functions 
using either fewer pairs of incomparable elements, or with a smaller number of
strict edges appearing.
So by repeated application of the above relations we may express any
$K_{(P, \omega)}^d(\xvec)$ as a linear combination of such generating functions
where every poset is totally ordered (a labeled chain) and
the labeling is natural. We use this later in Section~\ref{sec:nlwc}.

Similarly, one could also express any $K_{(P, \omega)}^d(\xvec)$
as a linear combination of labeled chains (not necessarily natural)
where all weights are $1$.
Such a $K_{(P, \omega)}(\xvec)$ for a labeled chain $(P, \omega)$
is the same as some \emph{fundamental quasisymmetric function} $\gessel_\alpha(\xvec)$; see \cite[Cor. 7.19.5]{StanleyEC2}.
\end{remark}

As a side note, $K^d_{(P, \omega)}(\xvec)$ for weighted labeled chains give rise to the 
\emph{universal quasisymmetric functions} defined in \cite[Section 3.2]{GrinbergVassilieva2021}.
These specialize to several of the standard families of quasisymmetric functions.

\subsection{Quasisymmetric power sums}

For a composition $\alpha$, let $\defin{z_\alpha} \coloneqq \prod_{i\geq1} i^{m_i} m_i!$,
where $m_i$ denotes the number of parts of $\alpha$ that are equal to $i$.
Given compositions $\alpha \leq \beta$ (denoting refinement)
so that $\beta = \beta_1 + \beta_2 + \dotsb + \beta_k$ 
and $\beta_i = \alpha_{i,1} + \alpha_{i,2} + \dotsb + \alpha_{i,\ell_i}$
for $1 \leq i \leq k$, let
\begin{equation}\label{eq:pidef}
\defin{\pi(\alpha,\beta)}
\coloneqq
\prod_{i=1}^{k} 
\alpha_{i,1} \cdot \left(\alpha_{i,1}+\alpha_{i,2}\right)\dotsm 
\left(\alpha_{i,1}+\alpha_{i,2} + \dotsb + \alpha_{i,\ell_i}  \right).
\end{equation}

The \defin{quasisymmetric power sum} $\Psi_\alpha$ (introduced and studied in \cite{BallantineDaughertyHicksMason2020},
but implicitly introduced
in \cite{GelfandKrobLascouxLeclercRetakhThibon1995} in connection with non-commutative symmetric functions) is then defined as
\begin{equation}\label{eq:PsiintoMonomial}
\defin{\Psi_\alpha(\xvec)} \coloneqq z_\alpha \sum_{\beta \geq \alpha } \frac{1}{\pi(\alpha,\beta)} M_\beta(\xvec),
\end{equation}
using the monomial quasisymmetric functions defined earlier in \eqref{eq:qmonomDef}.
The quasisymmetric power sums refine the power sum symmetric functions; 
\begin{equation}\label{eq:power sum_Psi_expansion}
p_\lambda(\xvec) = \sum_{\alpha \sim \lambda} \Psi_\alpha(\xvec)
\,,
\end{equation}
where the sum ranges over all compositions $\alpha$ whose parts re-arrange to $\lambda$;
see \cite[Prop.~3.19]{BallantineDaughertyHicksMason2020}.
We shall use the shorthand $\defin{\psih_\alpha} \coloneqq \frac{\Psi_\alpha}{z_\alpha}$
throughout the paper---these are the \defin{unnormalized quasisymmetric power sum functions}.

\subsection{\texorpdfstring{$(P,w)$}{(P,w)}-partitions and power sums}

We shall need a few definitions in order to state the next theorem.
Let $P$ be a poset on $n$ elements. Given a composition $\alpha \vDash n$,
let $\defin{\opsurj_{\alpha}(P)}$ 
denote the set of \defin{order-preserving surjections}
$f:P\to [\length(\alpha)]$ such 
that $|f^{-1}(i)| = \alpha_i$ for $i=1,\dotsc,\length(\alpha)$.
Note that $f^{-1}(i)$ is a sub-poset of $P$.
The set of all order-preserving surjections is denoted $\defin{\opsurj(P)}$
and we let $\defin{\wt(f)}$ be the composition $\alpha$ where $\alpha_i = |f^{-1}(i)|$.

\begin{definition}\label{def:pointed}
Let $P$ be a poset. An order-preserving surjection $f \in \opsurj_{\alpha}(P)$
is said to be \defin{pointed} if every sub-poset $f^{-1}(i)$ has a unique minimum.
We let $\defin{\opsurj^*(P)} \subseteq \opsurj(P)$ be the set of pointed order-preserving surjections.
\end{definition}

In \cite{AlexanderssonSulzgruber2019}, the following (after restating the result 
slightly using the notation in this article) results were proved.
Observe that the poset $(P,\omega)$ is assumed to be naturally labeled.
\begin{theorem}[{See \cite[Thm. 5.6]{AlexanderssonSulzgruber2019}}]\label{thm:ASFormula}
Let $(P, \omega)$ be a \emph{naturally labeled} poset.
Then
\begin{equation}
K_{(P, \omega)}(\xvec) = \sum_{\alpha} |\opsurj_\alpha^*(P)| \cdot  \psih_{\alpha}(\xvec) = \sum_{f \in \opsurj^*(P)} 
\psih_{\wt(f)}(\xvec).
\end{equation} 
\end{theorem}

\bigskip

There is a weighted generalization of the previous theorem.
Let $d$ be the weight function on $P$ as above,
and suppose $f \in \opsurj^*(P)$.
We let $\defin{\wt(f,d)}$ be the composition $(\alpha_1,\alpha_2,\dotsc,\alpha_\ell)$
where
\[
\alpha_i = \sum_{\poseta \in f^{-1}(i)} d(\poseta),
\]
that is, we total the weights of all vertices in the sub-poset $f^{-1}(i)$. We then have the following generalization of Theorem~\ref{thm:ASFormula}.

\begin{theorem}[{See \cite[Thm. 6.3]{AlexanderssonSulzgruber2019}}]\label{thm:ASFormula2}
Let $(P, \omega)$ be a naturally labeled poset and $d:P \to \setP$ be a weight function on $P$. Then
\begin{equation}
K_{(P, \omega)}^d(\xvec) = \sum_{\alpha} 
	\psih_{\alpha}(\xvec)
	\sum_{\substack{f \in \opsurj^*(P) \\ \wt(f,d) = \alpha}}
	\prod_{j=1}^{\ell(\alpha)}
	d\left( \min f^{-1}(j) \right),
	\end{equation} 
	where $f^{-1}(j)\subseteq P$ is a sub-poset.
\end{theorem}

Another result is due to Liu and Weselcouch~\cite{LiuWeselcouch2020}, where the authors give the 
quasisymmetric power sum expansion of $K_{(P, \omega)}(\xvec)$ for an \emph{arbitrary} labeling $\omega$.
In order to state their result, we need to introduce some additional terminology.

\begin{definition}[Generalized border strip]\label{def:gbs}
The labeled poset $(P, \omega)$ is a \defin{generalized border strip}
if there does not exist a chain $a<b<c$ in $P$ such that $\omega(a)<\omega(b)>\omega(c)$.
\end{definition}

\begin{remark}\label{rem:borderStrip}
Informally, a labeled poset is a generalized border strip if we cannot find 
two edges belonging in a chain of three elements, such that the 
bottom edge is natural, and the top edge is strict.
In particular, it is enough to consider the edges in the Hasse diagram of $P$.
\end{remark}

\begin{example}\label{gbstrip}
For example, the first labeled poset $(P, \omega)$ below is a generalized border strip,
while the second poset $(P', \omega)$ is not.
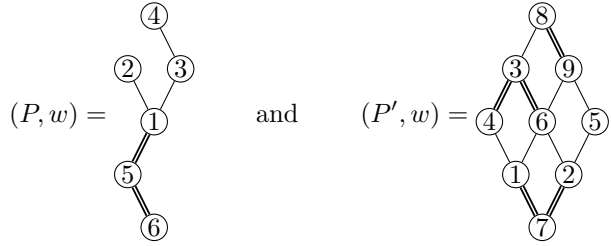
\[
(P, \omega) =
\begin{tikzpicture}
\draw
(1,1)node[circled](v1){6}
(1,5)node[circled](v4){1}
(0,7)node[circled](v6){2}
(0,3)node[circled](v7){5}
(2,7)node[circled](v8){3}
(1,9)node[circled](v9){4}
;
\draw (v4)--(v8);
\draw (v4)--(v6);
\draw (v8)--(v9);
\draw[strictEdge] (v1)--(v7);
\draw[strictEdge] (v7)--(v4);
\end{tikzpicture}
\qquad
\text{ and }
\qquad
	(P', \omega) = 
	\begin{tikzpicture}
	\draw
	(1,1)node[circled](v1){7}
	(2,3)node[circled](v3){2}
	(1,5)node[circled](v4){6}
	(3,5)node[circled](v5){5}
	(0,7)node[circled](v6){3}
	(0,3)node[circled](v7){1}
	(2,7)node[circled](v8){9}
	(1,9)node[circled](v9){8}
	(-1,5)node[circled](v10){4}
	;
	\draw (v7)--(v1)--(v3)--(v4)--(v6)--(v9);
	\draw (v3)--(v5);
	\draw (v4)--(v7);
	\draw (v5)--(v8);
	\draw (v4)--(v8);
	\draw (v7)--(v10);
	\draw[strictEdge] (v1)--(v3);
	\draw[strictEdge] (v4)--(v6);
	\draw[strictEdge] (v1)--(v7);
	\draw[strictEdge] (v8)--(v9);
	\draw[strictEdge] (v10)--(v6);
	\end{tikzpicture}
	\]
\end{example}

\begin{definition}[Liu--Weselcouch tagging]\label{LWdefn}

Now suppose $(P, \omega)$ is a generalized border strip. 
We define the \defin{LW-tagging} $\defin{\tagg}:P\to\{-1, 1^*, +1\}$ as follows.
If $x\in P$ is the bottom element of a strict edge, let $\tagg(x)=-1$.
If $x$ is the top element of a natural edge, let $\tagg(x)=1$.
Otherwise, let $\tagg(x)=1^*$. Note that $\tagg$ is well-defined (i.e.,~$x$ cannot be
the bottom element of a strict edge and the top element of a natural edge)
because $(P, \omega)$ is a generalized border strip. 
\end{definition}

\begin{definition}[Rooted generalized border strips]
We say $(P, \omega)$ is \defin{rooted} if it is a generalized 
border strip and there exists a unique $x\in P$ such that $\tagg(x)=1^*$.
Let 
\[
\defin{\sgn(P)} \coloneqq (-1)^{|\{x\in P:\tagg(x)=-1\}|}.
\]
We let $\defin{\rot(P)}$ denote the 
unique element tagged with $1^*$, whenever $(P, \omega)$ is rooted.
Note that $\sgn(P)$ and $\rot(P)$ also depend on the labeling $\omega$.
Observe that if $(P, \omega)$ is a naturally labeled poset, then its sign is $+1$,
and its root is its unique minimum if such an element exists.
\end{definition}

\begin{definition}[Border strip]
We almost exclusively consider \emph{rooted generalized border strips}, 
so we simply use the shorter \defin{border strip} for this term.
\end{definition}

\begin{example}
Each vertex in the generalized border strip 
given in Example~\ref{gbstrip} has the Liu--Weselcouch tagging as shown below:
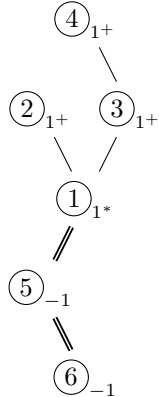
\[
\begin{tikzpicture}[scale=1.7]
\draw
(1,1)node(v1){\posetNode{6}{-1}}
(1,5)node(v4){\posetNode{1}{1^*}}
(0,7)node(v6){\posetNode{2}{+1}}
(0,3)node(v7){\posetNode{5}{-1}}
(2,7)node(v8){\posetNode{3}{+1}}
(1,9)node(v9){\posetNode{4}{+1}}
;
\draw (v4)--(v8);
\draw (v4)--(v6);
\draw (v8)--(v9);
\draw[strictEdge] (v1)--(v7);
\draw[strictEdge] (v7)--(v4);
\end{tikzpicture}
\]
Moreover, it is also rooted with $\defin{\rot(P)}$ being the element labeled $1$ whose LW-tagging is $1^{*}$.
\end{example}

\begin{definition}
Let $(P, \omega)$ be a labeled poset.
An order-preserving surjection $f \in \opsurj_{\alpha}(P)$
is said to be \defin{rooted} if every 
sub-poset $f^{-1}(i)$ is rooted (with respect to the induced labeling by $\omega$). The set of \defin{rooted order-preserving surjections} is denoted $\defin{\opsurj^*(P, \omega)}$.

It is easy to see that if $(P, \omega)$ is a naturally labeled poset,
then $f \in \opsurj^*(P, \omega)$
if and only if $f^{-1}(i)$ has a unique minimum,
so this notation agrees with Definition~\ref{def:pointed}.
\end{definition}

For general labelings, we have the following:
\begin{theorem}[{Liu--Weselcouch, \cite[Thm. 6.9]{LiuWeselcouch2020}}]\label{thm:liu}
Let $(P, \omega)$ be a labeled poset on $n$ vertices. 
Then
\begin{equation}
  K_{(P, \omega)}(\xvec) = \sum_{\alpha} 
  \psih_{\alpha}(\xvec)
  \sum_{\substack{f \in \opsurj^*(P, \omega) \\ \wt(f)=\alpha}} \prod_{i=1}^{\length(\alpha)} \sgn(f^{-1}(i)).
\end{equation}
\end{theorem}
In particular, the coefficient of $\psih_{(n)}$ in $K_{(P, \omega)}(\xvec)$
is $0$ unless $(P, \omega)$ is rooted, in which case the coefficient is 
$\pm 1$ according to the LW-tagging of $(P, \omega)$.

\subsection{The classical Murnaghan--Nakayama rule}\label{sec:MurnaghanNakayamaConnection}

In this subsection, we recall from \cite{LiuWeselcouch2020}
how the classical Murnaghan--Nakayama rule for Schur functions is obtained from Theorem~\ref{thm:liu}.
\medskip

The Schur function $\schur_{\lambda}(\xvec)$ is a symmetric function associated
with a partition $\lambda$; see \cite{StanleyEC2,Macdonald1995}.
It is easy to see that $\schur_{\lambda}(\xvec)$ is equal to the generating 
function $K_{(P, \omega)}(\xvec)$ if we choose an appropriate $(P, \omega)$. 
For example, the partition $\lambda = (6, 4, 3, 2, 2)$ corresponds 
to the labeled poset below:
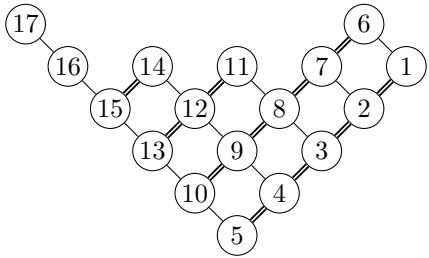
\[
\begin{tikzpicture}[scale=0.8]
\draw
(0,0)node[circled,minimum size=1.5em](v1){5}
(-2,2)node[circled,minimum size=1.5em](v2){10}
(2,2)node[circled,minimum size=1.5em](v3){4}
(0,4)node[circled,minimum size=1.5em](v4){9}
(-4,4)node[circled,minimum size=1.5em](v5){13}
(4,4)node[circled,minimum size=1.5em](v6){3}
(-2,6)node[circled,minimum size=1.5em](v7){12}
(2,6)node[circled,minimum size=1.5em](v8){8}
(0,8)node[circled,minimum size=1.5em](v9){11}
(6,6)node[circled,minimum size=1.5em](v10){2}
(-4,8)node[circled,minimum size=1.5em](v11){14}
(-6,6)node[circled,minimum size=1.5em](v12){15}
(-8,8)node[circled,minimum size=1.5em](v13){16}
(-10,10)node[circled,minimum size=1.5em](v14){17}
(4,8)node[circled,minimum size=1.5em](v15){7}
(8,8)node[circled,minimum size=1.5em](v16){1}
(6,10)node[circled,minimum size=1.5em](v17){6}
;
\draw (v1)--(v2)--(v5)--(v12)--(v13);
\draw (v3)--(v4)--(v7)--(v11);
\draw (v6)--(v8)--(v9);
\draw (v10)--(v15);
\draw (v16)--(v17);
\draw[strictEdge] (v1)--(v3)--(v6)--(v10)--(v16);
\draw[strictEdge] (v2)--(v4)--(v8)--(v15)--(v17);
\draw[strictEdge] (v5)--(v7)--(v9);
\draw[strictEdge] (v12)--(v11);
\draw (v13)--(v14);
\end{tikzpicture}
\]

The \defin{Murnaghan--Nakayama rule} (\cite{Murnaghan1937,Nakayama1940}) gives a
combinatorial formula for computing the coefficient $\chi_{\lambda}(\mu)$ in the power sum
expansion of Schur functions:
\begin{equation}\label{eq:mnrule}
\schur_{\lambda}(\xvec) = \sum_{\mu} \chi_{\lambda}(\mu) \frac{\powerSum_{\mu}(\xvec)}{z_\mu}.
\end{equation}

From \eqref{eq:power sum_Psi_expansion}, we recall that power
sum symmetric functions refine into the quasisymmetric power sum basis as
$\powerSum_{\lambda}(\xvec) = \sum_{\alpha \sim \lambda} \Psi_{\alpha}(\xvec)
$, where the sum runs over all compositions $\alpha$ that rearrange to $\lambda$.
Plugging this into \eqref{eq:mnrule}, we have that
\[
	\schur_{\lambda}(\xvec) =
	\sum_{\mu} \chi_{\lambda}(\mu) \cdot \frac{1}{z_{\mu}} \sum_{\alpha \sim \mu} \Psi_{\alpha}(\xvec) = \sum_{\alpha} \chi_{\lambda}(\alpha) \underbrace{\frac{\Psi_{\alpha}(\xvec)}{z_\alpha}}_{=\psih_{\alpha}(\xvec)}.
\]
Note that $\chi_{\lambda}(\alpha)=\chi_{\lambda}(\mu)$ if $\alpha \sim \mu$.
Recall now that the Murnaghan--Nakayama rule states that
\[
 \chi_{\lambda}(\mu) = \sum_{T \in \mathrm{BST}(\lambda,\mu)} (-1)^{\mathrm{height}(T)}
\]
where $\mathrm{BST}(\lambda,\mu)$ is the set of \defin{border strip tableaux} of shape $\lambda$
and type $\mu$, and $\mathrm{height}(T)$ is a certain statistic on these tableaux.
It is fairly straightforward to see that every $T \in \mathrm{BST}(\lambda,\mu)$
corresponds to some $f \in \opsurj^*(P, \omega)$ with $\wt(f)=\mu$, and
$\mathrm{height}(T)$ is precisely the total number of strict edges in all the border strips
determined by $f$. That is, Theorem~\ref{thm:liu} is indeed
a generalization of the classical Murnaghan--Nakayama rule; 
for more details, see \cite{LiuWeselcouch2020}.

\section{A weighted generalization}\label{sec:mainresult}

In this section, we unify the two formulas in Theorem~\ref{thm:ASFormula2} and Theorem~\ref{thm:liu}.
We present the proof of our main theorem below as follows: we begin by stating the
quasisymmetric power sum expansion of $K^d_{(P, \omega)}(x)$, and then show that this expansion satisfies the two recurrence relations in Lemma~\ref{lem:addEdge} and Lemma~\ref{lem:splitWeight}.

The argument proceeds inductively, reducing to the case of naturally labeled weighted chains.
This final step (see Lemma~\ref{lem:probabilistic}), is proved by a probabilistic interpretation.
\bigskip

The main theorem we want to prove is the following.
\begin{theorem}[Main theorem]\label{thm:main}
Let $(P, \omega, d)$ be a weighted labeled poset. Then
\begin{equation}\label{eq:mainTheorem}
  K^d_{(P, \omega)}(\xvec) = 
  \sum_{\alpha} 
  \psih_{\alpha}(\xvec)
  \sum_{\substack{f \in \opsurj^*(P, \omega) \\ \wt(f,d)=\alpha}}
  \prod_{i=1}^{\length(\alpha)} \sgn(f^{-1}(i)) d\left( \rot(f^{-1}(i)) \right).
\end{equation}
\end{theorem}
\begin{example}
Consider a $(P, \omega, d)$-partition from Example~\ref{gbstrip}, 
where each vertex has been assigned a weight:
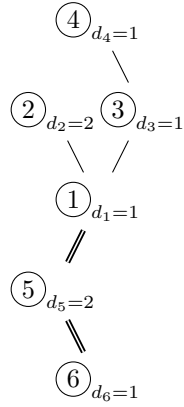
\[
\begin{tikzpicture}[xscale=2.2,yscale=1.6]
\draw
(1,1)node(v1){\posetNode{6}{d(\cdot)=1}}
(1,5)node(v4){\posetNode{1}{d(\cdot)=1}}
(0,7)node(v6){\posetNode{2}{d(\cdot)=2}}
(0,3)node(v7){\posetNode{5}{d(\cdot)=2}}
(2,7)node(v8){\posetNode{3}{d(\cdot)=1}}
(1,9)node(v9){\posetNode{4}{d(\cdot)=1}}
;
\draw (v4)--(v8);
\draw (v4)--(v6);
\draw (v8)--(v9);
\draw[strictEdge] (v1)--(v7);
\draw[strictEdge] (v7)--(v4);
\end{tikzpicture}
\]
The full $\psi$-expansion of $K^d_{(P, \omega)}(\xvec)$ is given by 
\begin{align*}
 &\psih_8-\psih_{17}-2 \psih_{35}+3 \psih_{62}+\psih_{71}+2 \psih_{125}-3 \psih_{152}-\psih_{161}-6 \psih_{332}-2 \psih_{341}+ \\
 &4 \psih_{422}+2 \psih_{512}+2 \psih_{521}+\psih_{611}+6 \psih_{1232}+2 \psih_{1241}-4 \psih_{1322}-2 \psih_{1412} +\\
 &-2 \psih_{1421}-\psih_{1511}-8 \psih_{3122}-4 \psih_{3212}-4 \psih_{3221}-2 \psih_{3311}+2 \psih_{4112}+2 \psih_{4121}+\\
 &2 \psih_{4211}+8 \psih_{12122}+4 \psih_{12212}+4 \psih_{12221}+2 \psih_{12311}-2 \psih_{13112}-2 \psih_{13121}-2 \psih_{13211}+\\
 &-4 \psih_{31112}-4 \psih_{31121}-4 \psih_{31211}+4 \psih_{121112}+4 \psih_{121121}+4 \psih_{121211}.
\end{align*}
Since the full labeled poset is a border strip, we have a non-zero coefficient of $\psih_8$.
As another example, the term $-3 \psih_{152}$ arises from two rooted order-preserving surjections,
one of which is 
\[
f^{-1}(1)=\{6^*\}, \quad f^{-1}(2)=\{1^*, 5, 2\}, \quad f^{-1}(3)=\{3^*, 4\}.
\]
Here, we refer to elements in the poset by their labels, and the roots are indicated by $*$.
These three border strips have signs $+1$, $-1$ and $+1$, respectively.
The product of weights of the roots is $1\cdot 1 \cdot 1$.

The second rooted order-preserving surjection is
\[
f^{-1}(1)=\{6^*\}, \quad f^{-1}(2)=\{1^*, 5, 3, 4\}, \quad f^{-1}(3)=\{2^*\},
\]
where the signs again are $+1$, $-1$ and $+1$, respectively.
The product of weights of the roots is $1\cdot 1 \cdot 2$.
In total, we get $(-1) +(-2) = -3$, which gives our coefficient.

The last three terms in the expansion above arise from 
the linear extensions of the poset.
\end{example}

\begin{remark}\label{rem:strict}
Note that the order-preserving surjections do not need
to adhere to the strict edges. In the first case above, we had $f^{-1}(2)=\{1^*, 5, 2\}$
so $f(1)=f(2)=f(5)=2$, even though there is a strict edge from $5$ to $1$
in the poset.
\end{remark}

The core of the proof of \eqref{eq:mainTheorem} is showing that the formula satisfies the two relations in 
Lemma~\ref{lem:addEdge} and Lemma~\ref{lem:splitWeight}.

\begin{proposition}\label{prop:addEdge}
The formula in \eqref{eq:mainTheorem} satisfies the recursion in Lemma~\ref{lem:addEdge}.
\end{proposition}
\begin{proof}
Consider two incomparable elements $a$ and $b$ 
in $P$---illustrated in the left-hand side in \eqref{eq:incomparable}---with weights $d(a)$ and $d(b)$, respectively.

The goal is to describe a
sign- and weight-preserving bijection
between rooted order-preserving surjections, which proves the following identity:
\begin{equation}\label{eq:incomparable}
\begin{tikzpicture}[scale=1.5]
\draw
(-2,0) node(extremeBL){}
( 4,0) node(extremeBR){}
(-2,6) node(extremeTL){}
( 4,6) node(extremeTR){}
(0,2) node(v1){\posetNode{a}{d(a)}}
(2,2) node(v2){\posetNode{b}{d(b)}}
(0,6)node(topLeft){ }
(1,6)node(topMid){ }
(2,6)node(topRight){ }
(0,0)node(botLeft){ }
(1,0)node(botMid){ }
(2,0)node(botRight){ };
\draw[gray,thick] (botLeft)--(v1);
\draw[gray,thick] (botMid)--(v1);
\draw[gray,thick] (botRight)--(v1);
\draw[gray,thick] (topLeft)--(v1);
\draw[gray,thick] (topMid)--(v1);
\draw[gray,thick] (topRight)--(v1);
\draw[gray,thick] (v2)--(topLeft);
\draw[gray,thick] (v2)--(topMid);
\draw[gray,thick] (v2)--(topRight);
\draw[gray,thick] (v2)--(botLeft);
\draw[gray,thick] (v2)--(botMid);
\draw[gray,thick] (v2)--(botRight);
\draw[dashed,blue] (v1)--(extremeBL);
\draw[dashed,blue] (v1)--(extremeTL);
\draw[dashed,red] (v2)--(extremeBR);
\draw[dashed,red] (v2)--(extremeTR);
\draw (1,-0.5) node(lbl) {$\opsurj^*(P, \omega)$};
\end{tikzpicture}
=
\begin{tikzpicture}[scale=1.5]
\draw
(-2,0) node(extremeBL){}
( 4,0) node(extremeBR){}
(-2,6) node(extremeTL){}
( 4,6) node(extremeTR){}
(0,2) node(v1){\posetNode{a}{d(a)}}
(2,4) node(v2){\posetNode{b}{d(b)}}
(0,6)node(topLeft){ }
(1,6)node(topMid){ }
(2,6)node(topRight){ }
(0,0)node(botLeft){ }
(1,0)node(botMid){ }
(2,0)node(botRight){ };
\draw[naturalEdge] (v1)--(v2);
\draw[gray,thick] (botLeft)--(v1);
\draw[gray,thick] (botMid)--(v1);
\draw[gray,thick] (botRight)--(v1);
\draw[gray,thick] (v2)--(topLeft);
\draw[gray,thick] (v2)--(topMid);
\draw[gray,thick] (v2)--(topRight);
\draw[dashed,blue] (v1)--(extremeBL);
\draw[dashed,blue] (v1)--(extremeTL);
\draw[dashed,red] (v2)--(extremeBR);
\draw[dashed,red] (v2)--(extremeTR);
\draw (1,-0.5) node(lbl) {$\opsurj^*(P', \omega)$};
\end{tikzpicture}
\!\!
+
\!\!
\begin{tikzpicture}[scale=1.5]
\draw
(-2,0) node(extremeBL){}
( 4,0) node(extremeBR){}
(-2,6) node(extremeTL){}
( 4,6) node(extremeTR){}
(0,4) node(v1){\posetNode{a}{d(a)}}
(2,2) node(v2){\posetNode{b}{d(b)}}
(0,6)node(topLeft){ }
(1,6)node(topMid){ }
(2,6)node(topRight){ }
(0,0)node(botLeft){ }
(1,0)node(botMid){ }
(2,0)node(botRight){ };
\draw[strictEdge] (v1)--(v2);
\draw[gray,thick] (botLeft)--(v2);
\draw[gray,thick] (botMid)--(v2);
\draw[gray,thick] (botRight)--(v2);
\draw[gray,thick] (v1)--(topLeft);
\draw[gray,thick] (v1)--(topMid);
\draw[gray,thick] (v1)--(topRight);
\draw[dashed,blue] (v1)--(extremeBL);
\draw[dashed,blue] (v1)--(extremeTL);
\draw[dashed,red] (v2)--(extremeBR);
\draw[dashed,red] (v2)--(extremeTR);
\draw (1,-0.5) node(lbl) {$\opsurj^*(P'', \omega)$};
\end{tikzpicture}
\end{equation}
\emph{
In what follows, we may sometimes ignore elements in $P$ which 
are only related to at most one of $a$ and $b$ --- i.e.~the dashed edges in the above figure, 
as these relations do affect the arguments.}

\bigskip 

Let $P$, $P'$ and $P''$ be as in Lemma~\ref{lem:addEdge} where $P'$ and $P''$ are the two ways of adding a relation between $a$ and $b$,
as illustrated in \eqref{eq:incomparable}. 
Let $f$ be an element of $\opsurj^*(P, \omega)$.
There are two cases to consider:
\begin{enumerate}[(i)]
 \item $a \in f^{-1}(i)$ and $b \in f^{-1}(j)$ for $i \neq j$, so $a$ and $b$
 belong to different border strips.
 In this case, $f$ is identically equal to either some $f' \in \opsurj^*(P', \omega)$ if $i<j$ or
 some $f'' \in \opsurj^*(P'', \omega)$ otherwise. This is because we add an edge between
 two \emph{different} border strips, so no tags within border strips are affected.
 Hence, we are done with this case.

 \item $a, b \in f^{-1}(i)$, and it suffices to prove the relation
 in the case when we only have a single border strip, i.e., when $(P, \omega)$ itself is a border strip.
\end{enumerate}
\emph{So it remains to treat case (ii), where we from now on
assume that $f \in \opsurj^*(P, \omega)$ is a single border strip.
We remark that in \cite{LiuWeselcouch2020},
the authors reduce to this case by a Hopf-algebraic argument.}
\medskip

We first show that there is an injection from objects
in the left-hand side to objects in the right-hand side.
Suppose we have a border strip for $(P, \omega)$ and that $a,b \in P$ are incomparable.
By symmetry, we may assume $\omega(a)< \omega(b)$.

There are a few signed cases to consider for an
element in $\opsurj^*(P, \omega)$,
namely the combinations of $\tagg(a),\tagg(b) \in \{-1, 1, 1^*\}$.
Note that not all combinations can appear.

We illustrate the various cases below---there might be
more edges present of course as sketched in \eqref{eq:incomparable}.
\begin{itemize}
 \item[\textbf{Case I:}] $\tagg(a)\in \{1, 1^*\}$ and $\tagg(b)=1$.
Here, we add a natural edge from $a$ to $b$ and obtain an element in $\opsurj^*(P',\omega)$.
Observe that if $\tagg(a)=1^*$, then there cannot be a natural edge to it from
some element below it in the poset, and that this map
preserves all tags. We illustrate this below:
\[
\begin{tikzpicture}[scale=1.7]
\draw
(0,0) node(v1){\posetNode{a}{d(a)}}
(2,0) node(v2){\posetNode{b}{d(b)}};
\draw[dashed] (v1)--(0,-2);
\draw[naturalEdge] (v2)--(2,-2);
\end{tikzpicture}
\quad
\longrightarrow
\quad
\begin{tikzpicture}[scale=1.7]
\draw
(0,0) node(v1){\posetNode{a}{d(a)}}
(2,2) node(v2){\posetNode{b}{d(b)}};
\draw[dashed] (v1)--(0,-2);
\draw[naturalEdge] (v2)--(2,-2);
\draw[naturalEdge] (v1)--(v2);
\end{tikzpicture}
\in \opsurj^*(P', \omega).
\]

\item[\textbf{Case II:}] $\tagg(a)\in \{1, 1^*\}$ and $\tagg(b)=-1$.
We insert a strict edge from $b$ to $a$ in this case:
 \[
\begin{tikzpicture}[scale=1.7]
\draw
(0,0) node(v1){\posetNode{a}{d(a)}}
(2,0) node(v2){\posetNode{b}{d(b)}};
\draw[dashed] (v1)--(0,-2);
\draw[strictEdge] (v2)--(2,2);
\end{tikzpicture}
\quad
\longrightarrow
\quad
\begin{tikzpicture}[scale=1.7]
\draw
(0,2) node(v1){\posetNode{a}{d(a)}}
(2,0) node(v2){\posetNode{b}{d(b)}};
\draw[dashed] (v1)--(0,-2);
\draw[strictEdge] (v2)--(2,2);
\draw[strictEdge] (v1)--(v2);
\end{tikzpicture}
\in \opsurj^*(P'', \omega).
\]

\item[\textbf{Case III:}] $\tagg(a)=-1$ and $\tagg(b)=1$.
In this situation, we insert a natural edge:
\[
\begin{tikzpicture}[scale=1.7]
\draw
(0,0) node(v1){\posetNode{a}{d(a)}}
(2,0) node(v2){\posetNode{b}{d(b)}};
\draw[strictEdge] (v1)--(0,2);
\draw[naturalEdge] (v2)--(2,-2);
\end{tikzpicture}
\quad
\longrightarrow
\quad
\begin{tikzpicture}[scale=1.7]
\draw
(0,-1) node(v1){\posetNode{a}{d(a)}}
(2,1) node(v2){\posetNode{b}{d(b)}};
\draw[strictEdge] (v1)--(0,2);
\draw[naturalEdge] (v2)--(2,-2);
\draw[naturalEdge] (v1)--(v2);
\end{tikzpicture}
\in \opsurj^*(P', \omega).
\]

\item[\textbf{Case IV:}] $\tagg(a)=-1$ and $\tagg(b)=-1$.
We insert a strict edge:
\[
\begin{tikzpicture}[scale=1.7]
\draw
(0,0) node(v1){\posetNode{a}{d(a)}}
(2,0) node(v2){\posetNode{b}{d(b)}};
\draw[strictEdge] (v1)--(0,2);
\draw[strictEdge] (v2)--(2,2);
\end{tikzpicture}
\quad
\longrightarrow
\quad
\begin{tikzpicture}[scale=1.7]
\draw
(0,0) node(v1){\posetNode{a}{d(a)}}
(2,-2) node(v2){\posetNode{b}{d(b)}};
\draw[strictEdge] (v1)--(0,2);
\draw[strictEdge] (v2)--(2,2);
\draw[strictEdge] (v1)--(v2);
\end{tikzpicture}
\in \opsurj^*(P'', \omega).
\]

\item[\textbf{Non-case:}] $\tagg(a)\in \{-1, 1\}$ and $\tagg(b)=1^*$.
This situation is not possible: Since $a$
is not the root, it is either below (in the Hasse diagram)
a vertex with smaller label,
or above some vertex with smaller label. Call this vertex
$c_1$. The same argument now applies to $c_1$,
so we get some vertex $c_2$ and so on.
Thus, there must be a sequence of vertices,
\[
 a, c_1, c_2,\dotsc,c_\ell = b
\]
ending at the unique root $b$.
But then $\omega(a)>\omega(c_1)>\dotsb > \omega(b)$,
contradicting $\omega(a)<\omega(b)$.
\end{itemize}
In conclusion, we have showed that \emph{if} $f \in \opsurj^*(P, \omega)$ is a border strip then (by adding an edge) it can be mapped
to a border strip in \emph{exactly one of} $\opsurj^*(P', \omega)$, $\opsurj^*(P'', \omega)$,
in such a way that all tags and weights are preserved.

\bigskip 

\emph{But we are not done yet!}
The cases in the right-hand side which are \emph{not} in
the image of the injection are the ones where
either $(P', \omega)$ or $(P'', \omega)$ is a border strip,
but removal of the edge $(a, b)$ results in a generalized
border strip with two roots. 
Let us consider such border strips:
\begin{itemize}
 \item Suppose we have $f \in \opsurj^*(P', \omega)$ 
 (where $a$ is the root and $(a,b)$ is a weak edge) 
 where removal of $(a, b)$ results in a generalized border strip with two roots.
 Only the tagging at vertices $a$ and $b$ can change, so both $a$ and $b$ must be 
 roots after removing this edge. 
 
	This can only happen if all edges above $a$ are weak and all edges below $a$ 
	are strict --- the same must hold for $b$, see Figure~\ref{fig:2roots}.
 
 \item Similarly, consider $f \in \opsurj^*(P'', \omega)$ 
 (where $a$ is the root and $(a,b)$ is a strict edge) where again removal of $(a,b)$
 results in two roots.
 We have the exact same conditions on the labeling as in previous case.
\end{itemize}
This shows that there is a bijection between the above two cases, see Figure~\ref{fig:2roots}.
\begin{figure}[!ht]
\begin{align*}
\substack{
\begin{tikzpicture}[scale=1.5]
\draw
(0,-1) node(v1){\posetNode{a}{d(a)}}
(2,1) node(v2){\posetNode{b}{d(b)}};
\draw[blue] (v1)--(0,3);
\draw[blue,strictEdge] (v1)--(0,-3);
\draw[blue] (v2)--(2,3);
\draw[blue,strictEdge] (v2)--(2,-3);
\draw[naturalEdge] (v1)--(v2);
\end{tikzpicture}
\\
\in \opsurj^*(P', \omega)
}
\longleftrightarrow
\begin{tikzpicture}[scale=1.5]
\draw
(0,-1) node(v1){\posetNode{a}{d(a)}}
(2,1) node(v2){\posetNode{b}{d(b)}};
\draw[blue] (v1)--(0,3);
\draw[blue,strictEdge] (v1)--(0,-3);
\draw[blue] (v2)--(2,3);
\draw[blue,strictEdge] (v2)--(2,-3);
\end{tikzpicture}
\longleftrightarrow
\begin{tikzpicture}[scale=1.5]
\draw
(0,1) node(v1){\posetNode{a}{d(a)}}
(2,-1) node(v2){\posetNode{b}{d(b)}};
\draw[blue] (v1)--(0,3);
\draw[blue,strictEdge] (v1)--(0,-3);
\draw[blue] (v2)--(2,3);
\draw[blue,strictEdge] (v2)--(2,-3);
\end{tikzpicture}
\longleftrightarrow
\substack{
\begin{tikzpicture}[scale=1.5]
\draw
(0,1) node(v1){\posetNode{a}{d(a)}}
(2,-1) node(v2){\posetNode{b}{d(b)}};
\draw[blue] (v1)--(0,3);
\draw[blue,strictEdge] (v1)--(0,-3);
\draw[blue] (v2)--(2,3);
\draw[blue,strictEdge] (v2)--(2,-3);
\draw[strictEdge] (v1)--(v2);
\end{tikzpicture}
\\
\in \opsurj^*(P'', \omega)}
\end{align*}
\caption{The step-by-step bijection 
from $\opsurj^*(P', \omega)$ (left) to $\opsurj^*(P'', \omega)$ (right) 
for the cases where removal of $(a,b)$ results in two roots (the two middle posets).
All edges except $(a,b)$ are optional and do not have any effect, 
but the type (strict/weak) must be as indicated.
Note, $\tagg(b)=1$ in the leftmost poset, but $\tagg(b)=-1$ in the rightmost.
}\label{fig:2roots}
\end{figure}

Since $\tagg(b)$ is $1$ in $P'$ and $-1$ in $P''$, the border strips have opposite signs.
We therefore have a sign-reversing involution on all elements not hit by the previous injection.
\end{proof}

\subsection{Splitting the weights}

\begin{proposition}\label{prop:splitWeight}
The formula in \eqref{eq:mainTheorem} satisfies the identity in Lemma~\ref{lem:splitWeight}.
\end{proposition}
\begin{proof}
Recall that the identity we want to verify is illustrated as follows:
\begin{equation}\label{eq:splitWeightProof}
 \begin{tikzpicture}[scale=1.7]
\draw
(1,3)node(v1){\posetNode{i}{d_1+d_2}}
(0,6)node(topLeft){ }
(1,6)node(topMid){ }
(2,6)node(topRight){ }
(0,0)node(botLeft){ }
(1,0)node(botMid){ }
(2,0)node(botRight){ };
\draw[thick,gray] (botLeft)--(v1)--(topLeft);
\draw[thick,gray] (botMid)--(v1)--(topMid);
\draw[thick,gray] (botRight)--(v1)--(topRight);
\draw (1,-0.5) node(lbl) {$\opsurj^*(P, \omega)$};
\end{tikzpicture}
=
 \begin{tikzpicture}[scale=1.7]
\draw
(0,2) node(v1){\posetNode{i_-}{d_1}}
(2,4) node(v2){\posetNode{i_+}{d_2}}
(0,6)node(topLeft){ }
(1,6)node(topMid){ }
(2,6)node(topRight){ }
(0,0)node(botLeft){ }
(1,0)node(botMid){ }
(2,0)node(botRight){ };
\draw[naturalEdge] (v1)--(v2);
\draw[thick,gray] (botLeft)--(v1);
\draw[thick,gray] (botMid)--(v1);
\draw[thick,gray] (botRight)--(v1);
\draw[thick,gray] (v2)--(topLeft);
\draw[thick,gray] (v2)--(topMid);
\draw[thick,gray] (v2)--(topRight);
\draw (1,-0.5) node(lbl) {$\opsurj^*(P', \omega')$};
\end{tikzpicture}
-
 \begin{tikzpicture}[scale=1.7]
\draw
(2,4) node(v1){\posetNode{i_-}{d_2}}
(0,2) node(v2){\posetNode{i_+}{d_1}}
(0,6)node(topLeft){ }
(1,6)node(topMid){ }
(2,6)node(topRight){ }
(0,0)node(botLeft){ }
(1,0)node(botMid){ }
(2,0)node(botRight){ };
\draw[strictEdge] (v1)--(v2);
\draw[thick,gray] (botLeft)--(v2);
\draw[thick,gray] (botMid)--(v2);
\draw[thick,gray] (botRight)--(v2);
\draw[thick,gray] (v1)--(topLeft);
\draw[thick,gray] (v1)--(topMid);
\draw[thick,gray] (v1)--(topRight);
\draw (1,-0.5) node(lbl) {$\opsurj^*(P', \omega'')$};
\end{tikzpicture}
\end{equation}
Recall from Lemma~\ref{lem:splitWeight}, that the labelings $\omega'$ and $\omega''$ satisfy
\[
 \omega'(a_-)=i_-, \;
 \omega'(a_+)=i_+, \quad  
 \omega''(a_-)=i_+, \; \omega''(a_+)=i_-.
\]

Suppose we have some $f \in \opsurj^*(P, \omega)$ for the poset in the left-hand side,
and where $a \in P$ is the element labeled $i$.
We consider the border strip $S$ containing $a$.
There are three cases on the left-hand side:
\begin{enumerate}[(i)]
\item\label{c:1} Vertex $a$ (labeled $i$) is not a root, 
and there is some natural edge connecting to $a$ in $S$ from below.
\item\label{c:2} Vertex $a$ is not a root, and there is a strict edge connecting to $a$ in $S$ from above.
\item\label{c:3} Vertex $a$ is the root.
 Top edges are natural and bottom edges are strict (if any of these are present).
\end{enumerate}
Note that cases (\ref{c:1}) and (\ref{c:2}) are disjoint since $S$ is a border strip.
\medskip

In case (\ref{c:1}), $f$ induces a border-strip $f' \in \opsurj^*(P', \omega')$
that agrees with $f$ outside $a$ and where we have $f'(a_-) = f'(a_+) = f(a)$.
This $f'$ appears in in the first term in the right-hand side of \eqref{eq:splitWeightProof},
and the LW-tags for $a_-$ and $a_+$ are both $+1$ for $f'$.

The surjection $f$ does not induce a border strip for the second term.
Note that if $b\neq a$ is the root in the left-hand side,
then $b$ is still a root for the induced border strip in the right-hand side
so the weights of the roots in both sides of \eqref{eq:splitWeightProof} are the same.

In case (\ref{c:2}), $f$ does not induce a border strip for $\opsurj^*(P', \omega')$,
but $f$ does induce a border strip for the second term.
The LW-tag of $a$ is $-1$ and in the right-hand side,
we have that both $a_-$ and $a_+$ have LW-tag $-1$.
Thus the total sign on the right-hand side term is the same as for $f$.
As in previous case, the weights of the roots agree.

Finally, in case (\ref{c:3}), $f$ now contributes with weight $(d_1+d_2)$
on the left-hand side.
In the right-hand side, the induced border strip for the first term has $a_-$ (labeled $i_-$) as the root,
contributing with $d_1$. Moreover, the induced border strip for the 
second term has $a_+$ (labeled $i_-$) as the root and one extra element tagged $-1$.
Hence, this term contributes with weight $d_2$.
The LW-tags of all other elements are preserved, so the signs match.

\bigskip 
So we now have a sign-preserving and weight-preserving \emph{injection}
\[
  \opsurj^*(P, \omega) \xhookrightarrow{}  \opsurj^*(P', \omega') \sqcup \opsurj^*(P', \omega''),
\]
and whose image covers all cases where $a_-$ and $a_+$  are in the same border strip.
However, in the right-hand side we need to consider the cases where 
the elements $a_-$ and $a_+$ belong to different border strips.
But these terms simply cancel since an edge between different border strips
has no impact on the sign or the weight.
\end{proof}

By Remark~\ref{rem:recursion}, Proposition~\ref{prop:addEdge}
and Proposition~\ref{prop:splitWeight} it suffices to show that
the main formula \eqref{eq:mainTheorem} holds in the case with labeled chains and no weights.
This special case is simply \cite[Thm. 3.1]{AlexanderssonSulzgruber2019}.
Alternatively, we can get rid of the weights using Proposition~\ref{prop:splitWeight}
and directly apply the unweighted rule in \cite{LiuWeselcouch2020}.
For the sake of completeness, we finish the proof by giving a new purely combinatorial 
argument which avoids the reference to Hopf algebras or fundamental quasisymmetric functions.

\subsection{Naturally labeled weighted chains}\label{sec:nlwc}

By Lemma~\ref{lem:addEdge} and Lemma~\ref{lem:splitWeight},
it follows that any $K^d_{(P, \omega)}(\xvec)$ can be expressed
as some linear combination $\sum_j K^{d^j}_{C_j}(\xvec)$
where each $C_j$ is a naturally labeled chain and each $d^j$ is a weight function $d^j : C_j \to \setN^+$.

Thus, in order to finish the proof of \eqref{eq:mainTheorem} (by an induction argument) 
it suffices to show that it holds for this base case.
Suppose $(P, \omega)$ is a chain with vertex 
weights $d = (d_1,d_2,\dotsc,d_\ell)$ (bottom to top) with total weight $n$.
Then 
\begin{equation}\label{eq:LHS1}
 K^d_{(P, \omega)}(\xvec) = \sum_{i_1 \leq i_2 \leq \dotsb \leq i_\ell} x_{i_1}^{d_1} \dotsm x_{i_\ell}^{d_\ell}
 =
  \sum_{\beta \geq d} M_\beta(\xvec).
\end{equation}
Since we are in the naturally labeled case, we need to show that this agrees with 
\begin{equation}\label{eq:weightedChainInPsi}
  \sum_{f \in \opsurj^*(P, \omega) }   \psih_{\wt(f,d)}(\xvec)
  \prod_{i} d\left( \rot(f^{-1}(i)) \right).
\end{equation}
Since $(P, \omega)$ is a chain,
the compositions $\wt(f, d)$ that appear in \eqref{eq:weightedChainInPsi}
correspond precisely to coarsenings of the composition $d = (d_1, d_2, \dots, d_\ell)$. That is, any $\psih_{\alpha}(\xvec)$ appearing in \eqref{eq:weightedChainInPsi} must be obtained from a coarsening of the form
\begin{equation}\label{eq:alphaCoarsening}
\alpha = 
  (
  \underbrace{d_1+d_2+\dotsb+d_{i_2-1}}_{\alpha_1},
  \underbrace{d_{i_2}+\dotsb + d_{i_3-1}}_{\alpha_2},\dotsc, 
  \underbrace{d_{i_k} +\dotsb + d_\ell}_{\alpha_k} ).
\end{equation}
Hence, \eqref{eq:weightedChainInPsi} is obtained by summing 
over all coarsenings $\alpha$ of $d$:
\[
  \sum_{\alpha \geq d} \psih_{\alpha}(\xvec) \cdot d_1\cdot d_{i_2}\cdot  \dotsm \cdot d_{i_k}
\]
where the $i$'s are as in \eqref{eq:alphaCoarsening}.
By the definition in \eqref{eq:PsiintoMonomial}, this now equals
\begin{equation}\label{eq:RHS1}
 \sum_{\alpha \geq d} \sum_{\beta \geq \alpha} \frac{M_{\beta}(\xvec)}{\pi(\alpha,\beta)}
\cdot d_1\cdot d_{i_2}\cdot  \dotsm \cdot d_{i_k}.
\end{equation}

We now compare the coefficients of $M_\beta(\xvec)$ in \eqref{eq:RHS1} and \eqref{eq:LHS1},
and we must show that for a fixed $\beta \geq d$, we have
\begin{equation}\label{eq:toBe1}
 1 = 
\sum_{\substack{\alpha \\
d \leq \alpha \leq \beta }} \frac{1}{\pi(\alpha,\beta)}
\cdot d_1\cdot d_{i_2}\cdot  \dotsm \cdot d_{i_k}.
\end{equation}
The $\alpha$'s are refinements of $\beta$ and we can refine each part of 
$\beta$ independently. Moreover, since $\pi(\alpha,\beta)$ is a product (see \eqref{eq:pidef})
over the parts of $\beta$, it suffices to show \eqref{eq:toBe1} for the one-part case $\beta = (n)$.
That is, we want to show
\[
 \sum_{\alpha \geq d} d_1\cdot d_{i_2}\cdot  \dotsm \cdot d_{i_k}
 \prod_{j=1}^k \left(\alpha_1 + \alpha_2 + \dotsb + \alpha_j \right)^{-1} = 1.
\]

It then suffices to prove the following lemma.
\begin{lemma}\label{lem:probabilistic}
Let $d=(d_1,d_2,\dotsc, d_\ell)$ be a composition. 
Then 
\begin{equation}\label{eq:omegaAlphaTerm}
  \sum_{\alpha \geq d} 
 \prod_{j=1}^k \frac{d_{i_j}}{\alpha_1 + \alpha_2 + \dotsb + \alpha_j} = 1,
\end{equation}
where $i_2,i_3,\dotsc,i_k$ are given by
\[
 \alpha = (d_1+d_2+\dotsb+d_{i_2-1},
  d_{i_2}+\dotsb + d_{i_3-1},\dotsc, 
d_{i_k} +\dotsb + d_\ell ), \text{where } i_1 = 1.
\]
\end{lemma}
\begin{proof}
We provide a probabilistic proof of this identity by defining an appropriate probability space, $\Omega$.
This space consists of all possible integer vectors:
	\[
	\defin{\Omega} \coloneqq \left\{ (a_1,a_2,\dotsc,a_\ell) : 1 \leq a_i \leq d_1+d_2+\dotsb +d_i \right\}.
	\]
	The total number of such vectors is  
	\[
	d_1 (d_1+d_2)\dotsm (d_1+d_2+\dotsb+d_\ell),
	\]
and we uniformly choose a vector from this set.
To visualize this space, we introduce an $\ell \times \ell$ staircase diagram, where the $j^\text{th}$ column contains $j$ boxes.
Each vector in $\Omega$ corresponds to selecting one box from each column.
The boxes in row $i$ (counting from the bottom) have weight $d_i$,
meaning that \defin{a box in row $i$ represents $d_i$ choices}.
We say that this is a \defin{staircase diagram},
and a selection of boxes then represents \emph{several} events in $\Omega$.
The number of events is simply the product of the weights of the selected boxes.

For instance, the $7\times 7$ staircase diagrams in
Figure~\ref{fig:staircaseDiagrams} illustrates two different outcomes.
\begin{figure}[!ht]
\begin{ytableau}
\none[d_7] & \none &\none & \none & \none &\none&\none&\none &\, \\
\none[d_6] & \none &\none & \none & \none &\none&\none&\,&\checkmark \\
\none[d_5] & \none &\none & \none & \none &\none&\,&\,&\, \\
\none[d_4] & \none &\none & \none & \none &\,&\,&\,&\, \\
\none[d_3] & \none &\none & \none & \checkmark &\,&\checkmark&\,&\, \\
\none[d_2] & \none &\none & \, & \, &\checkmark &\,&\checkmark &\, \\
\none[d_1] & \none &\checkmark & \checkmark & \, &\,&\,&\,&\, \\
\end{ytableau}
\hfill
\begin{ytableau}
\none[d_7] & \none &\none & \none & \none &\none&\none&\none &\, \\
\none[d_6] & \none &\none & \none & \none &\none&\none&\,&\, \\
\none[d_5] & \none &\none & \none & \none &\none&\,&\,&\; \\
\none[d_4] & \none &\none & \none & \none &\,&\,&\,&\, \\
\none[d_3] & \none &\none & \none & \checkmark &\,&\checkmark&\,&\, \\
\none[d_2] & \none &\none & \checkmark & \, &\checkmark &\,&\checkmark &\, \\
\none[d_1] & \none &\checkmark & \; & \, &\,&\,&\,&\checkmark \\
\end{ytableau}
\caption{Two examples of staircase diagrams illustrating sets of events in the probability space $\Omega$.
The first diagram represents $d_1^2\cdot d_2^2 \cdot d_3^2 \cdot d_6$
different integer vectors in $\Omega$.
}\label{fig:staircaseDiagrams}
\end{figure}

For each integer composition $\beta$ of $\ell$,
we have a corresponding coarsening $\alpha$ of $(d_1,d_2,\dotsc,d_\ell)$.
For example, 
\begin{equation}\label{eq:coarsening}
 \alpha = (d_1+d_2, d_3+d_4+d_5,d_6) \quad \longleftrightarrow \quad  \beta = (2,3,1).
\end{equation}
Our goal now is to partition the probability space $\Omega$ 
into disjoint events $\Omega_\beta$ indexed by integer compositions $\beta$ of $\ell$.

Given the staircase diagram of a vector in $\Omega$, we find which $\Omega_\beta$ it belongs to as follows.
\begin{itemize}
\item[(a)] We start with the staircase diagram and consider the selected box in the last column.
The position of the selected box from top will be the last part, $\beta_m$, of the composition $\beta$. 
	
\item[(b)] 
We remove the last $\beta_m$ columns from the staircase diagram,
and then proceed inductively, from step (a). 
This process gives the remaining parts $(\beta_1,\beta_2,\dotsc,\beta_{m-1})$.
\end{itemize}
By construction, it is clear that every staircase diagram belongs to some $\Omega_\beta$,
and that $\Omega$ is the disjoint union of the various $\Omega_\beta$.

The events belonging to $\Omega_\beta$ can be illustrated using a staircase diagram,
but we shade the columns that were not ``the last column'' in step (a) at any point, see Figure~\ref{fig:shading}.
For instance, the staircase diagram in Figure~\ref{fig:staircaseDiagrams} (left)
above belongs to $\Omega_{(2,3,2)}$.
\begin{figure}[!ht]
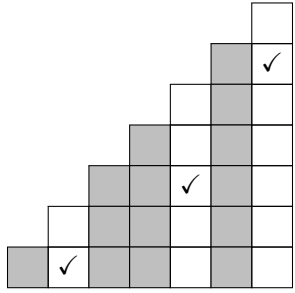

\begin{ytableau}
\none & \none & \none &\none&\none&\none &\, \\
\none & \none & \none &\none&\none &*(lightgray)&\checkmark \\
\none & \none & \none &\none&\,&*(lightgray)&\, \\
\none & \none & \none &*(lightgray)&\,&*(lightgray)&\, \\
\none & \none & *(lightgray) &*(lightgray)&\checkmark&*(lightgray)&\, \\
\none & \, & *(lightgray) &*(lightgray) &\,&*(lightgray) &\, \\
*(lightgray) & \checkmark & *(lightgray) &*(lightgray)&\,&*(lightgray)&\, \\
\end{ytableau}
\caption{The staircase diagrams belonging to $\Omega_{(2,3,2)}$ must have the marked boxes selected,
but we are free to chose any box in each shaded column.
}\label{fig:shading}
\end{figure}
The probability of a vector belonging to $\Omega_{(2,3,2)}$ is then
given by 
\[
\frac{d_1}{(d_1+d_2)} \cdot \frac{d_3}{(d_1+d_2+d_3+d_4+d_5)} \cdot \frac{d_6}{(d_1+d_2+d_3+d_4+d_5+d_6+d_7)}.
\]
This is precisely one of the terms in the left-hand side of \eqref{eq:omegaAlphaTerm},
and it is straightforward to see that the terms in  \eqref{eq:omegaAlphaTerm}
are precisely the probabilities for belonging to the $\Omega_\beta$s.
The desired result now follows since the total sum of probabilities is $1$.
\end{proof}

\begin{example}\label{ex:figg3}
Let us consider a case where $\ell=3$. 
Then the probability space consists of the set of all possible integer
vectors $\Omega = \{(a_1, a_2, a_3)\}$,
where  $1\leq a_1 \leq d_1$, $1\leq a_2 \leq d_1+d_2$ and $1\leq a_3 \leq d_1+d_2+d_3$.
The total number of such vectors is $d_1(d_1+d_2)(d_1+d_2+d_3)$.

Using $3\times 3$ staircase diagrams, each such vector corresponds to selecting one 
box from each column, where boxes in row $1$, row $2$ and row $3$ (from the bottom) have weights $d_1, d_2$ and $d_3$, respectively. 
Partitioning this space into disjoint events indexed by integer 
compositions $\beta$ of $3$ gives rise to the events belonging to $\Omega_\beta$ as illustrated in 
Figure~\ref{fig:fig3}.

\begin{figure}[!ht]
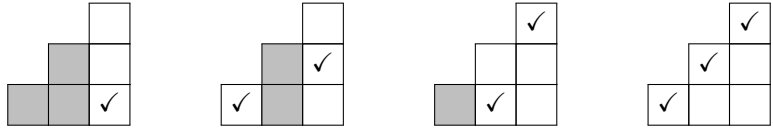

	\begin{ytableau}
	\none&\none &\, \\
	\none &*(lightgray)&\, \\
	*(lightgray)&*(lightgray)&\checkmark \\
	\end{ytableau}
\begin{ytableau}
\none &\none&\none&\none &\, \\
\none &\none&	\none &*(lightgray)&\checkmark \\
\none &\none&	\checkmark &*(lightgray)&\, \\
\end{ytableau}
\begin{ytableau}
	 \none &\none&\none&\none &\checkmark \\
	 \none &\none&\none &\,&\,\\
	 \none &\none&*(lightgray)&\checkmark &\, \\
\end{ytableau}
\begin{ytableau}
	 \none &\none&\none&\none &\checkmark \\
	 \none &\none&\none &\checkmark &\, \\
	 \none &\none&\checkmark &\,&\, \\
\end{ytableau}
	\caption{The staircase diagrams for $\Omega_{(3)}$, $\Omega_{(1, 2)}$, $\Omega_{(2, 1)}$, $\Omega_{(1, 1, 1)}$ respectively.}\label{fig:fig3}
\end{figure}

The probabilities of a vector belonging to each of the events listed in Figure~\ref{fig:fig3} are
\begin{align*}
&\frac{d_1}{(d_1+d_2+d_3)},\quad  \frac{d_1}{(d_1)} \cdot \frac{d_2}{(d_1+d_2+d_3)},\quad \frac{d_1}{(d_1+d_2)} \cdot \frac{d_3}{(d_1+d_2+d_3)},\\
& \frac{d_1}{(d_1)} \cdot \frac{d_2}{(d_1+d_2)}\cdot \frac{d_3}{(d_1+d_2+d_3)},
\end{align*}
respectively. It is easily verified that the total sum of these probabilities is $1$.
\end{example}

We can generalize the argument in Lemma~\ref{lem:probabilistic} 
to prove the following $q$-analog of \eqref{eq:omegaAlphaTerm}.
\begin{proposition}
Let $d=(d_1,d_2,\dotsc, d_\ell)$ be a composition. Then
\begin{equation}\label{eq:omegaAlphaTermQ}
\sum_{\alpha \geq d} 
 \prod_{j=1}^k \frac{q^{\alpha_1 + \alpha_2 + \dotsb + \alpha_{j-1}}  [d_{i_j}]_q  }{[\alpha_1 + \alpha_2 + \dotsb + \alpha_j]_q} = 1,
\end{equation}
where \defin{$[m]_q \coloneqq 1+q+q^2 + \dotsb + q^{m-1}$} and $i_j$
is as in Lemma~\ref{lem:probabilistic}.
\end{proposition}
\begin{proof}[Proof sketch:]
We let each outcome $(a_1,\dotsc,a_\ell)$ in $\Omega$ 
have total weight $q^{a_1+a_2+\dotsb +a_\ell - \ell}$ instead,
so that the weighted sum of all vectors is
\[
 [d_1]_q [d_1+d_2]_q \dotsm [d_1+d_2+\dotsb+d_\ell]_q.
\]
The first factor corresponds to the weighted $d_1$ choices represented by the box in the first column;
the second factor is the weighted choices in the second column and so on.
So in particular, in the last column there are in total $d_1+d_2+\dotsb+d_\ell$
way to choose $a_\ell$ and the total sum of the weights of these 
is $[d_1+d_2+\dotsb+d_\ell]_q$. 
A box in row $d_{i_j}$ has $\alpha_1+\alpha_2+\dotsb+\alpha_{j-1}$
boxes below, so the sum of the weights of the options represented by a box in row $d_{i_j}$
is precisely the numerator in \eqref{eq:omegaAlphaTermQ}.

The rest of the proof follows from the same argument as above.
\end{proof}

We also have a connection with linear extensions of certain trees.
We note that the denominators appearing in \eqref{eq:omegaAlphaTerm}
can be interpreted as \defin{hook products} of certain trees.
\begin{definition}
Given $(d_1,\dotsc,d_\ell)$ and a 
composition $\beta \vDash \ell$, we define the \defin{$\beta$-tree} as follows.
Let $\beta = (\beta_{1}, \beta_{2},\dotsc ,\beta_{k})$ for some $k \leq \ell$,
and let $\alpha$ be the coarsening of $(d_1,\dotsc,d_\ell)$ determined by $\beta$
as in \eqref{eq:coarsening}.

Consider the tree $T$ having $k$ internal vertices which 
appear at levels $1,2,\dotsc,k$ from the bottom.
To the $i^\thsup$ internal vertex, we then add in 
total $\alpha_i - 1$ leaves. We can partition the leaves of $v_1$ into blocks 
of sizes given by $d_1-1$, $d_2$, $d_3$ and so on, in this order;
see Figure~\ref{fig:fig4}.
The total number of vertices in $T$ is then $d_1+d_2+\dotsc+d_\ell$.
For instance, the $\beta$-tree corresponding to $\Omega_{(2, 3, 2)}$ is shown below.

\begin{figure}[!ht]
\scalebox{0.75}{
	\begin{tikzpicture}[level distance=2.5cm,line width=1pt,
	level 1/.style={sibling distance=2.5cm},
	level 2/.style={sibling distance=2.5cm},
	level 3/.style={sibling distance=2.5cm}, scale=0.8]
	\node[hookTreeNodeInternal] {$v_3$}
	child {
		node[hookTreeNodeInternal] {$v_2$}
		child {node[hookTreeNodeInternal] {$v_1$}
			child {node[hookTreeNode] {$d_{1}-1$}}
			child {node[hookTreeNode] {$d_{2}$}}
		}
		child {node[hookTreeNode] {$d_{3}-1$}}
		child {node[hookTreeNode] {$d_{4}$}}
		child {node[hookTreeNode] {$d_{5}$}}
	}
	child {
		node[hookTreeNode] {$d_{6}-1$}
	}
	child {
		node[hookTreeNode] {$d_{7}$}
	};
	\end{tikzpicture}
}
	\caption{The $\beta$-tree for $\Omega_{(2, 3, 2)}$.
	Each leaf in the figure represents a number of leafs in the actual $\beta$-tree.
	For example, the node $d_4$ in the illustration represents $d_4$ leaves of $v_2$. The hook values for the internal vertices $v_1, v_2, v_3$ are $d_1+d_2$,  $d_1+d_2+d_3+d_4+d_5$,
	and $d_1+d_2+\dotsb+d_7$, respectively.
	}\label{fig:fig4}
\end{figure}
\end{definition}

The \defin{hook value} of a vertex in a tree is defined as $1$ plus its total number 
of descendants. 
Let $\mathrm{LinearExtensions}(\beta,d)$ denote the set of linear extensions 
of the $\beta$-tree with weights $d$.
The hook formula for counting linear extensions of trees; 
see \cite[Chap. 5.1.4, Ex. 20]{Knuth1998ArtOfProgramming}, then tells us that 
\[
 |\mathrm{LinearExtensions}(\beta,d)|\cdot \prod_{j=1}^k \left(\alpha_1 + \alpha_2 + \dotsb + \alpha_j \right) = (d_1+d_2+\dotsb + d_\ell)!.
\]
Using this relation in \eqref{eq:omegaAlphaTerm}, we arrive at
\[
 \sum_{\beta \vDash \ell} |\mathrm{LinearExtensions}(\beta,d)| \cdot \prod_{j=1}^k d_{i_j} = (d_1+d_2+\dotsb + d_\ell)!.
\]
The above formula is remarkable in the sense that the left-hand 
side is a sum of terms involving linear extensions,
while the right-hand side is a simple product. 
\begin{problem}
Give a direct bijective proof of this identity. Is there a $q$-analog? Can it be generalized?
Note that when all $d_i$ are equal to $1$, then the above formula in essence
tells us that permutations can be expressed uniquely in
cycle form; see \cite[Proof of Thm.~2.6, Case I]{AlexanderssonSulzgruber2019}.
Briefly, internal vertices in the tree correspond to the smallest element in each cycle,
and leaves are the other elements in the cycle.
\end{problem}

\section{Open problems}\label{sec:problems}

\subsection{Comparison with the cylindrical Murnaghana--Nakayama rule}

The weak and strict inequalities in the generating function $K_{(P, \omega)}(\xvec)$
are determined by the labeling $\omega$.
However,  not all combinations of weak and strict edges in the poset 
can be realized using a labeling; see P.~McNamara~\cite{McNamara2006}. 
One can instead consider the more general notion of \emph{oriented posets} where arbitrary edges 
can be marked as strict. For example, the class of \emph{cylindrical Schur functions} 
are realized as the generating functions of colorings of such oriented posets.

There is a generalization of the Murnaghan--Nakayama rule for cylindrical Schur functions;
see \cite{AlexanderssonOguz2023x}. Notably, for this cylindrical rule, we require that 
some border strips come with weights (not just a sign).
This indicates that it is significantly harder to 
find a rule for the $\psi$-expansion of general oriented posets.

\begin{problem}
Find a Murnaghan--Nakayama rule that gives the $\psi$-expansion of the generating function of 
oriented posets.
This would unify both the $\psi$-expansion of $K_{(P, \omega)}(\xvec)$ and 
the Murnaghan--Nakayama rule for cylindrical Schur functions.
\end{problem}

\subsection{Rectangular compositions and cancellation-free sums}

The classical Murnaghan--Nakayama rule is cancellation-free for Schur functions
whenever we want to compute the coefficient of $\powerSum_{k^\ell}$.
This result goes back to James--Kerber~\cite{JamesKerber1984}.
Their method generalizes to skew Schur functions without much effort; see~\cite{AlexanderssonPfannererRubeyUhlin2020}.

It is therefore natural to ask if the Liu--Weselcouch formula~\eqref{thm:liu}
is cancellation-free when all border strips have the same size.
That is, if $\alpha = (k, k,\dotsc,k)$, does the product in
\[
	\sum_{\substack{f \in \opsurj^*(P, \omega) \\ \wt(f)=\alpha}} \prod_{i=1}^{\length(\alpha)} \sgn(f^{-1}(i))
\]
always have the same sign, independently of $f$?
The answer to this question is in general no, as the next example shows.

\begin{example}[Cancellation]
Consider the labeled poset $(P, \omega)$ in Figure~\ref{fig:fig5},
where we consider the coefficient of $\psih_{22}$.
We need to partition $P$ into two rooted border strips of size $2$.
There are two ways to accomplish this.
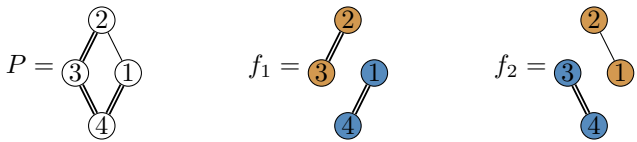
\begin{figure}[!ht]
\begin{equation*}
P=
\begin{tikzpicture}
\draw
(1,5)node[circled](v4){4}
(0,7)node[circled](v6){3}
(2,7)node[circled](v8){1}
(1,9)node[circled](v9){2}
;

\draw (v8)--(v9);
\draw[strictEdge] (v4)--(v6);
\draw[strictEdge] (v4)--(v8);
\draw[strictEdge] (v6)--(v9);
\end{tikzpicture}
\qquad
\qquad
f_1=
	\begin{tikzpicture}
	\draw
	(1,5)node[circled,fill=pLightBlue](v4){4}
	(0,7)node[circled,fill=pLightOrange](v6){3}
	(2,7)node[circled,fill=pLightBlue](v8){1}
	(1,9)node[circled,fill=pLightOrange](v9){2}
	;
	
	\draw[strictEdge] (v4)--(v8);
	\draw[strictEdge] (v6)--(v9);
	\end{tikzpicture}
\qquad
\qquad
f_2=
	\begin{tikzpicture}
	\draw
	(1,5)node[circled,fill=pLightBlue](v4){4}
	(0,7)node[circled,fill=pLightBlue](v6){3}
	(2,7)node[circled,fill=pLightOrange](v8){1}
	(1,9)node[circled,fill=pLightOrange](v9){2};
	\draw[strictEdge] (v4)--(v6);
	\draw (v9)--(v8);
	\end{tikzpicture}
\end{equation*}
\caption{The poset $P$ and the two different rooted order-preserving surjections of $P$.}\label{fig:fig5}
\end{figure}
Note that $\sgn(f_1) = (-1)^2 = 1$ while $\sgn(f_2) = (-1)\cdot 1 = -1$,
so the coefficient of $\psih_{22}$ is $0$ in the expansion of $K_{(P, \omega)}(\xvec)$.
\end{example}

\begin{problem}
 Find a large natural family of labeled posets where the Murnaghan--Nakayama rule is cancellation-free
 when computing the coefficient of quasisymmetric power sums indexed by rectangular shapes.
\end{problem}

\begin{problem}
Are the \emph{set-valued $P$-partitions} (see e.g. \cite[Def. 5.4]{LamPylyavskyy2007combinatorialhopf}) positive in the quasisymmetric power sum basis?
Can the Liu--Weselcouch formula be extended to set-valued $(P, \omega)$-partitions?
For example, there is a type of Murnaghan--Nakayama rule for stable Grothendieck polynomials of Grassmann type; see \cite{NguyenHiepSonThuy2023}.
\end{problem}

\subsection*{Acknowledgements}

The second author acknowledges the International Science Programme (ISP) for their 
generous support through the Mikael Passare Postdoctoral Award, 
which funded her post-doctoral visit at Stockholm University's Department of Mathematics.

We are also grateful for the insightful comments from Jörgen Backelin.

Finally, we sincerely thank the anonymous reviewers for their
detailed and insightful comments,
which greatly improved the clarity and quality of this paper.

\bibliographystyle{amsalpha}
\bibliography{bibliography.bib}

\end{document}